\documentclass[a4paper,11pt]{article}
\usepackage[latin1]{inputenc}
\usepackage[english]{babel}
\usepackage{amssymb}
\usepackage{amsmath}
\usepackage{amsthm}
\usepackage{mathrsfs}
\usepackage[all]{xy}
\usepackage[dvips]{graphicx}
\usepackage{url}
\title{About the stability of the tangent bundle of $\mathbb{P}^n$ restricted to a surface}
\author{Chiara Camere \thanks{Laboratoire J.-A. Dieudonné, UMR CNRS 6621, Université de Nice-Sophia Antipolis, Parc Valrose, 06108 Nice Cedex 02 - \textit{E-mail:} \texttt{camere56@unice.fr}}}
\begin{document}
\date{}
\maketitle

\newtheorem{defi}{Definition}
\newtheorem{lem}{Lemma}
\newtheorem{pro}{Proposition}
\newtheorem{teo}{Theorem}
\newtheorem{cor}{Corollary}
\newtheorem{rmk}{Remark}
\newcommand{\of}{\mathcal{O}}
\newcommand{\ff}{\mathscr{F}}
\newcommand{\pn}{\mathbb{P}}
\newcommand{\rk}{\mathrm{rk}\,}
\begin{abstract}
Let $ X $ be a smooth projective surface over $ \mathbb{C} $ and let $ L $ be a line bundle on $ X $ generated by its global sections. Let $ \phi _L:X\longrightarrow\pn ^r $ be the morphism associated to L; we investigate the $ \mu- $stability of $ \phi _L^*T_{\pn ^r} $ with respect to $ L $ when $ X $ is either a regular surface with $ p_g=0 $, a K3 surface or an abelian surface. In particular, we show that $ \phi _L^*T_{\pn ^r} $ is $ \mu- $stable when $ X $ is K3 and $ L $ is ample and when $ X $ is abelian and $ L^2\geq 14 $.
\end{abstract}
\section{Introduction}

Given a line bundle $L$ generated by its global sections on a smooth projective variety $X$, one can consider the kernel of the evaluation map
\begin{equation}\label{eq:syzX}
\xymatrix{
0\ar@{->}[r]&M_L\ar@{->}[r]&H^0(X,L)\otimes\mathcal{O}_X \ar@{->}[r]&L\ar@{->}[r]&0}
\end{equation}
and its dual $E_L=M_L^*$.

The stability of this bundle is equivalent to that of $\phi _L^*T_{\pn ^r}$, where $ \phi _L:X\longrightarrow\pn ^r $ is the morphism associated to L. It has been studied in the case of a curve by Paranjape in \cite{PR} with Ramanan and in his Ph.D. thesis \cite{P}; in particular, the latter contains the statements on which rely all our results contained in a former paper \cite{C} and in this one. Later Ein and Lazarsfeld showed in \cite{EL} that $M_L$ is stable if $\deg L>2g$ and Beauville investigated the case of degree $2g$ in \cite{Be1}.

The aim of this paper is to study this problem in the case of projective surfaces. Here we consider the $\mu-$stability of a sheaf with respect to a chosen linear series $H$, which generalises the definition given in the case of curves: a vector bundle $E$ is said to be $\mu-$stable with respect to $H$ if for each proper quotient sheaf $F$ we have $\mu(F)>\mu(E)$, where $\mu(F)= \frac{c_1(F).H^{n-1}}{\rk F}$ is the slope of $F$ (see \cite{huy}).

After studying these vector bundles in Section \ref{sim}, we gather some results which hold on curves in Section \ref{s1} and then in Section \ref{s3} we obtain some results about regular surfaces, including the following 
\begin{teo} \label{p:must}
Let $X$ be a smooth projective K3 surface over $\mathbb{C} $ and let $ L$ be an ample line bundle generated by its global sections on $ X $; then the vector bundle $E_L$ is $\mu-$stable with respect to $L$.
\end{teo}

Finally, in Section \ref{s6} we study the case of abelian surfaces, showing the following
\begin{teo}\label{abel}
Let $X$ be a smooth projective abelian surface over $\mathbb{C}$ and let $L$ be a line bundle on $X$ generated by its global sections such that $L^2\geq 14$. Then the vector bundle $E_L$ is $\mu-$stable with respect to $L$.
\end{teo}
\section{Simplicity and rigidity of $E_L$}\label{sim}

Let us briefly recall the geometric interpretation of $ E_L  $: since $L$ is generated by its global sections, the morphism $ \phi _L:X\longrightarrow\pn (H^0(L))\simeq\pn ^r $ is well-defined and we have $L=\phi _L^*\of _{\pn ^r}(1) $; thus, from the dual sequence of  (\ref{eq:syzX}) and from the well-known Euler exact sequence
it follows that $ E_L=\phi _L^*T_{\pn ^r}\otimes L^* $ and the stability of $E_L $ is equivalent to the stability of $\phi _L^*T_{\pn ^r}$.

In the next sections we will deal with the problem of whether or not these bundles are $\mu$-stable, but let us first of all underline that they satisfy in almost any case a less strong property, the simplicity.

\begin{pro}
Let $X$ be a smooth projective variety and $L$ be a big line bundle generated by its global sections on $ X $; if $\dim X\geq 2$ then $E_L$ is simple. 
 \end{pro}
\textbf{Proof.} If we tensor with $E_L$ the short exact sequence (\ref{eq:syzX}) in cohomology we get
\begin{equation}\label{seqcoml}
\xymatrix{
0\ar@{->}[r]&H^0(M_L\otimes E_L)\ar@{->}[r]&H^0(L)\otimes H^0(E_L) \ar@{->}[r]^{\alpha}&H^0(L\otimes E_L)\ar@{->}[r]&\\
\ar@{->}[r]&H^1(M_L\otimes E_L)\ar@{->}[r]&H^0(L)\otimes H^1(E_L) \ar@{->}[r]&\cdots}
\end{equation} 
Since $H^0(L^*)\cong H^1(L^*)\cong 0$ by Ramanujam-Kodaira vanishing theorem (see \cite{Mu}), we also have $H^0(L)^*\cong H^0(E_L)$. Now, by tensoring the dual sequence of (\ref{eq:syzX}) with $L$ we obtain in cohomology
\begin{equation}\label{simple}
\xymatrix{
0\ar@{->}[r]&H^0(\of _X)\ar@{->}[r]&H^0(L)\otimes H^0(L)^* \ar@{->}[r]^{\alpha}&H^0(L\otimes E_L)\ar@{->}[r]&H^1(\of _X)\ar@{->}[r]&\cdots}
\end{equation}
where the morphism $\alpha$ is the same morphism as in (\ref{seqcoml}).
Hence $H^0(M_L\otimes E_L)\cong H^0(\of _X)\cong \mathbb{C}$, i.e. $E_L$ is simple.\qed
\\

In the case of regular surfaces, under mild assumptions, which hold for example if $X$ is a K3 surface, they are also rigid, hence providing an example of an exceptional vector bundle on such a surface.
\begin{pro}
Let $X$ be a smooth projective regular surface and $L$ as above; if the multiplication map $ H^0(K_X)\otimes H^0(L)\rightarrow H^0(K_X\otimes L)$ is surjective, then $E_L$ is rigid. 
\end{pro}

\textbf{Proof.}  The morphism $\alpha$ in sequence (\ref{simple}) is surjective because $X$ is regular. Let us show that $H^1(E_L)\cong 0$: indeed, by tensoring (\ref{eq:syzX}) with $K_X$ in cohomology we get 
\[
\xymatrix{
0\ar@{->}[r]&H^0(M_L\otimes K_X)\ar@{->}[r]&H^0(L)\otimes H^0(K_X) \ar@{->}[r]^{\varphi}&H^0(L\otimes K_X)\ar@{->}[r]&\\
\ar@{->}[r]&H^1(M_L\otimes K_X)\ar@{->}[r]&H^0(L)\otimes H^1(K_X)=0 }
\]
Since we assumed $\varphi$ surjective, we have $H^1(E_L)\cong H^1(M_L\otimes K_X)\cong 0$ by the duality theorem. Then from the exact sequence (\ref{seqcoml}) it follows that $\mathrm{Ext}^1(E_L,E_L)\cong H^1(M_L\otimes E_L)\cong 0$, i.e. $E_L$ is rigid.\qed 

\section{Some results on vector bundles on curves}\label{s1}
Let us briefly recall some facts about vector bundles on curves.
In a former paper \cite{C} we showed the following
\begin{teo}\label{p:pa}
Let $C$ be a smooth projective curve of genus $g\geq 2$ over an algebraically closed field $k$ and let $L$ be a line bundle on $C$ generated by its global sections such that $\deg L\geq 2g-c(C) $. Then:
\begin{enumerate}
 \item $ E_L $ is semi-stable;
 \item $ E_L $ is stable except when $\deg L= 2g$ and either $C$ is hyperelliptic or $L\cong K(p+q)$ with $p,q\in C$.
\end{enumerate}
\end{teo}

In the case $L=K_C$ more was already known: in \cite{PR} Paranjape and Ramanan showed the following
\begin{teo}\label{p:para}
Let $C$ be a smooth projective curve of genus $g\geq 2$ over $\mathbb{C}$; $E_{K_C}$ is always semistable and it is also stable if $C$ is not hyperelliptic.
\end{teo}

The proof of Theorem \ref{p:pa} was essentially based on the following lemma, shown by Paranjape in \cite{P}.
\begin{lem}\label{p:p32}
Let $F$ be a vector bundle on $C$ generated by its global sections and such that $ H^0(C,F^*)=0 $; then 
$ \deg F\geq \rk F+g-h^1(C,\det F) $.
Moreover, if $ h^1(C,\det F)\geq 2 $ then $ \deg F\geq 2\rk F+c(\det F)\geq 2\rk F+c(C) $.
\end{lem}

\section{About regular surfaces}\label{s3}
Before restricting to the case of regular surfaces, let us see a few statements which hold for every surface.

\begin{lem}\label{zp}
Let $F$ be a vector bundle of rank 2 generated by its global sections on a smooth projective surface $X$ and assume moreover that $h^0(\det F)=2$. Then there is a short exact sequence
\begin{equation}\label{sk31}
 \xymatrix{
0\ar@{->}[r]&\of _X\ar@{->}[r]^s&F\ar@{->}[r]&\det F\ar@{->}[r]&0}\end{equation}
\end{lem}
\textbf{Proof.} We cannot have $F=\of _X^2$ because $h^0(\det F)=2$; then, since $ F $ is of rank 2 generated by its global sections, we can suppose $ h^0(F)\geq 3 $. Then there is a section $ s\in H^0(X,F) $ which is zero only in a finite number of  points and we have the following short exact sequence \begin{equation}\label{sk3}
\xymatrix{
0\ar@{->}[r]&\of _X\ar@{->}[r]^s&F\ar@{->}[r]&\mathcal{I}_Z\det F\ar@{->}[r]&0}\end{equation}
where $ Z $ is the zero locus of $ s $.
In cohomology we obtain
\[\xymatrix{
0\ar@{->}[r]&H^0(X,\of _X)\ar@{->}[r]&H^0(X,F)\ar@{->}[r]&H^0(X,\mathcal{I}_Z\det F)\ar@{->}[r]&\cdots}\]
Since $ h^0(F)\geq 3 $, we get $ h^0(\mathcal{I}_Z\det F) \geq 2$, but $ h^0(\mathcal{I}_Z\det F) \leq h^0(\det F)=2$.  Since $\det F$ is generated by its global sections, from $ h^0(\mathcal{I}_Z\det F) = h^0(\det F)=2 $ it follows that $ \mathcal{I}_Z\det F= \det F $ and $ Z=\varnothing $. Therefore the sequence (\ref{sk3}) becomes (\ref{sk31}).
\qed

\begin{pro}\label{pdiag}
Let $X$ be a smooth projective surface over $\mathbb{C}$ and let $L$ be a line bundle on $X$ generated by its global sections. Let $C$ be a smooth irreducible curve on $X$ such that $H^1(L\otimes\of _X(-C))=0$. Then $(E_L)_{|C}=E_{(L_{|C})}\oplus\of _C ^{r}$, with  $r=h^0(L\otimes\of _X(-C))$.
\end{pro}
\textbf{Proof.} 
Tensoring the exact sequence
\[ 
\xymatrix{
0\ar@{->}[r]&\of _X(-C)\ar@{->}[r]&\mathcal{O}_X \ar@{->}[r]&\of _C\ar@{->}[r]&0}\]
with $L$, we get
\[ 
\xymatrix{
0\ar@{->}[r]&L\otimes\of _X(-C)\ar@{->}[r]&L \ar@{->}[r]&L _{|C}\ar@{->}[r]&0}\]
and hence in cohomology we have
$$
\xymatrix{
0\ar@{->}[r]&H^0(X,L\otimes\of _X(-C))\ar@{->}[r]&H^0(X,L)\ar@{->}[r]&H^0(X,L_{|C})\ar@{->}[r]&0}
$$

So we have the following diagram
\begin{equation}\label{supdia}
 \xymatrix{
&0\ar@{->}[d]&0\ar@{->}[d]&0\ar@{->}[d]&\\
0\ar@{->}[r]&L_{|C}^* \ar@{->}[r]\ar@{->}[d]&H^0(X,L_{|C})^*\otimes\mathcal{O}_C\ar@{->}[r]\ar@{->}[d]&E_{(L_{|C})}\ar@{->}[r]\ar@{->}[d]&0\\
0\ar@{->}[r]&L_{|C}^* \ar@{->}[r]\ar@{->}[d]&H^0(X,L)^*\otimes\mathcal{O}_C \ar@{->}[r]^-{e_L}\ar@{->}[d]& (E_L)_{|C} \ar@{->}[r]\ar@{->}[d]&0\\
0\ar@{->}[r]& 0\ar@{->}[r]\ar@{->}[d]&\of _C ^{r}\ar@{->}[r]\ar@{->}[d]&\of _C ^{r}\ar@{->}[r]\ar@{->}[d]&0\\
&0&0&0&}
\end{equation}
By the snake lemma, the third column is exact. Moreover, the sequence splits and $(E_L)_{|C}=E_{(L_{|C})}\oplus\of _C ^{r}$.\qed\\

\begin{cor}\label{regpg0}
Let $X$ be a smooth projective regular surface over $\mathbb{C}$ such that $p_g=0$ and let $C$ be a smooth irreducible curve on $X$ of genus $g\geq 2$ such that $L=\of _X(K_X+C)$ is generated by its global sections; then $E_L$ is $\mu-$semistable with respect to $C$ and it is also stable if $c(C)>0$.
\end{cor}
\textbf{Proof.} By Proposition \ref{pdiag} $(E_L)_{|C}\cong E_{(L_{|C})}$, since $r=p_g=0$; on the other hand, $L_{|C}=K_C$, so the statement follows from Theorem \ref{p:para}.\qed\\

When $r\neq 0$, the restriction to the curve is no longer semistable, but in the case of K3 surfaces this is enough to gain the $\mu-$stability.\\

\textbf{Proof of Theorem \ref{p:must}.} Let $C\in |L|$ be a smooth irreducible curve of genus $g\geq 2$. By Proposition \ref{pdiag} we have $(E_L)_{|C}=E_{K_C}\oplus\of _C $, since $L_{|C}\cong K_C$; moreover $\mu(E_L)=\frac{2g-2}{g}< 2$. Let us suppose that $g\geq 3$: if $g=2$ then $C$ is hyperelliptic and we will deal with the case $c(C)=0$ later. Let $F$ be a quotient sheaf of $E_L$ of rank $ 0<\rk F<g $; then $F_{|C}$ is a quotient of $(E_L)_{|C} $. There is a diagram of the form
$$
\xymatrix{0\ar@{->}[r]&\of_C\ar@{->}[r]\ar@{->}[d]&(E_L)_{|C} \ar@{->}[r]\ar@{->}[d]&E_{K_C}\ar@{->}[r]\ar@{->}[d]&0\\
0\ar@{->}[r]&W\ar@{->}[r]\ar@{->}[d]&F_{|C} \ar@{->}[r]\ar@{->}[d]&G\oplus\tau\ar@{->}[r]\ar@{->}[d]&0\\
&0&0&0&
}
$$
where $G$ is a vector bundle generated by its global sections, $W$ is either $\of _C$ or $0$ and $\tau$ is a torsion sheaf on $C$, hence $\deg W= 0$ and $\deg \tau\geq 0$. So we get $\mu (F)=\frac{\deg G+\deg\tau}{\rk F}$.

\begin{itemize}
\item If $ \rk G=0 $, then $\rk (F)=1$ and we always have $\mu (F)\geq 2$. Indeed, otherwise it would be $ F= \of_X(D)$ with $ D>0 $ an effective base-point free divisor such that $ D.C=0$ or $1 $; we cannot have $D.C=0$, since $D$ is nef, hence $D^2\geq 0$, but by the Hodge index theorem we would have $D^2<0$, which is a contradiction. If $D.C=1$, by the Hodge index theorem we get $D^2=0$, hence $D=kE$ with $k\geq 1$ and $E$ an elliptic curve; in fact, we have $k=1$ because $D.C=1$, so $h^0(D)=2$ and $|D|$ is a pencil; then, since $C.D=1$, $C$ would be a section and $C^2<0$, impossible.

\item If $ \rk G>0 $, then $G$ is generated by its global sections such that $H^0(C,G^*)=0$; the hypothesis of Lemma \ref{p:p32} then hold and, since $\mu (F)\geq \frac{\deg G}{\rk F}$, we have:
\begin{enumerate}
\item if $h^1(\det G)<2$, since $g\geq 3$, then $$ \mu (F)\geq 1+\frac{g-2}{\rk G+1}>1+\frac{g-2}{g}=\mu (E_L).$$

\item If $h^1(\det G)\geq 2$, then $$\mu (F)\geq 2+\frac{c(\det G)+\deg\tau-2}{\rk G+1}\geq 2>\mu (E_L)$$ if $c(\det G)
 \geq 2$, in particular if $c(C)\geq 2$, but also if $c(\det G)=1$ and $\deg\tau>0$.

\end{enumerate}
\end{itemize}
This shows that $\mu (F)>\mu (E_L)$ in the case $c(C)\geq 2$. \\

We now deal with the case $c(C)=1$. We can repeat the above proof by applying Lemma \ref{p:p32} and it does not work only if $h^1(\det G)\geq 2$, $\deg \tau=0$ and $c(\det G)=1$. If $ g=3 $ then $ \mu(E_L)=\frac{4}{3} $ and we always have $ \mu(F)> \frac{4}{3}$.

 From now on we assume $ g\geq 4 $; then either the curve is trigonal or a smooth plane quintic of genus $g=6$ (see \cite{Mar}).

\begin{enumerate}
\item If there is a $\mathfrak{g}^1 _3$ on $C$, the only line bundles which compute the Clifford index are $\of _C(\mathfrak{g}^1 _3)$ and $\of_C(K_C-\mathfrak{g}^1 _3)$. 
  \begin{enumerate}
\item If $\det G=\of _C(\mathfrak{g}^1 _3)$, since $h^1(\det G)\geq 2$, by Lemma \ref{p:p32} we have $\deg G\geq 2\rk G+1$, hence in this case $\rk G= 1$.  Then $\rk F=2$ and $\det F_{|C}=\of _C(\mathfrak{g}^1 _3)$; it follows that $\det F=\of _X(D)$ with $D.C=3$. By the Hodge index theorem then, since $g\geq 4$, we have $D^2\leq \frac{9}{2g-2}<2$, so $D^2=0$ and $D=kE$ with $k\geq 1$ and $E$ an elliptic curve; since $D.C=3$ and $C.E\geq 2$, this implies $k=1$ and  $ h^0(\of _X(D))=2 $; by Lemma \ref{zp}, it follows from $ h^1(\det F^*)=0=\mathrm{Ext}^1(\of_X,\det F) $ that $ F=\of_X\oplus \det F $, hence $ h^0(F^*)>0 $, which is impossible.
\item If $\det G=\of _C(K_C-\mathfrak{g}^1 _3)$ we have $\deg G=2g-5$ and $\rk G\leq g-3$ by Lemma \ref{p:p32}, hence $$ \mu (F)\geq \frac{2g-5}{\rk G+1}\geq\frac{2g-5}{g-2}=2-\frac{1}{g-2}>\mu (E_L)$$if $g>4$. If $g=4$ we have $\deg G=3$ and we fall in the former case.
\end{enumerate}
\item If there is a $\mathfrak{g}^2 _5$ on $C$, the genus is $g=6$ and the only line bundle which computes the Clifford index is $\of _C(\mathfrak{g}^2 _5)\cong\of  _C(K_C-\mathfrak{g}^2 _5)$. 

If $\det G=\of _C(\mathfrak{g}^2 _5)$, since $h^1(\det G)\geq 2$, by Lemma \ref{p:p32} $\deg G\geq 2\rk G+1$, hence $\rk G\leq 2$ and $\rk F\leq 3$. Therefore we get $$ \mu (F)=\frac{5}{\rk G+1}\geq\frac{5}{3}=\mu (E_L)$$ 
Let us investigate whether equality can hold or not; suppose that $\rk F=3$. Since $F$ is of rank $>2$ generated by its global sections, there is a short exact sequence
\begin{equation}\label{sf}
 \xymatrix{
0\ar@{->}[r]&\of _X\ar@{->}[r]&F\ar@{->}[r]&V\ar@{->}[r]&0}\end{equation}
with $V$ of rank 2 generated by its global sections such that $\det V=\det F=\of _X(D)$ with $D.C=5$. By the Hodge index theorem then  $ D^2\leq2 $;  however the case $D^2=2$ cannot occur, since otherwise $(C-2D)^2=-2$ and by Riemann-Roch theorem at least one between $C-2D$ and $2D-C$ would be effective, contradicting $(C-2D).C=0$ and the ampleness of $C$.
If $D^2=0$, then $D=kE$ with $k\geq 1$ and $E$ an elliptic curve; since $D.C=5$ and $C.E\geq 2$, this implies  $k=1$ and  $ h^0(\of _X(D))=2 $, so by Lemma \ref{zp} there is a short exact sequence
\[\xymatrix{
0\ar@{->}[r]&\of _X\ar@{->}[r]^s&V\ar@{->}[r]&\det V\ar@{->}[r]&0}\]
and in cohomology we obtain $h^1(V^*)=h^1(V)=0$.
As a consequence we have $\mathrm{Ext}^1(\of_X,V)=0$ and $F=\of_X\oplus V$, impossible since it would imply $h^0(F^*)>0$.
\end{enumerate}

Then $\mu (F)>\mu (E_L)$ also if $c(C)=1$.\\

Suppose now that $C$ is a hyperelliptic curve; in this case (see \cite{Be}, pag.129), the morphism $\phi _L:X\longrightarrow \pn ^g$ induces a double covering  $\pi: X\longrightarrow F$ where $F\subset\pn ^g$ is a rational surface of degree $g-1$ which is either smooth or a cone over a rational normal curve. If $g=2$ then $F=\pn ^2$  (see \cite{Be}, pag.129) and it is well-known that its tangent bundle is $\mu-$stable (see \cite{huy} Section 1.4) with respect to $\of _{\pn ^2}(1)$. If $g\geq 3$, let $i:F\hookrightarrow\pn ^g$ be the embedding and $H=i^*\of _{\pn ^g}(1)$ the ample hyperplane section of $F$ such that $\pi^*H=L$; we have $H^2=g-1$.

On the surface $F$ we have the short exact sequence 
\begin{equation}\label{mh}
 \xymatrix{
0\ar@{->}[r]&H^*\ar@{->}[r]&H^0(F,H)^*\otimes\mathcal{O}_F \ar@{->}[r]& E_H\ar@{->}[r]&0}\end{equation}
We know that the curve $H$ is rational, so $p_a(H)=0$; we consider a smooth curve $\Gamma\in |2H|$. By the adjunction formula we have $0=p_a(H)=1+\frac{1}{2}(H^2+H.K_F)$, so we get $H.K_F=-H^2-2=-g-1$; using the adjunction formula once more we then obtain
\[p_a(\Gamma)=1+\frac{1}{2}(\Gamma ^2+\Gamma .K_F)=1+2H^2+H.K_F=g-2\]

Since $g\geq 3$ we have $p_a(\Gamma)\geq 1$. Since $H$ is ample, we deduce $H^0(F,\of _F(-H))=H^1(F,\of _F(-H))=0$ (see \cite{Mu}).
Then from the short exact sequence
\[ 
\xymatrix{
0\ar@{->}[r]&\of _F(H-\Gamma)\ar@{->}[r]&\of _F(H) \ar@{->}[r]&\of _{\Gamma}(H)\ar@{->}[r]&0}\]
and from the associated cohomology sequence it follows that $H^0(F,\of _F(H))\cong H^0(F,\of _{\Gamma}(H))$, hence $(E_H)_{|\Gamma}=E_{\of _{\Gamma}(H)}$.

Moreover, $\deg \of _{\Gamma}(H)=H.\Gamma=2g-2>2p_a(\Gamma)=2g-4$. Since $\of _{\Gamma}(H)$ is a line bundle on a smooth projective curve $\Gamma$ of genus $\geq 1$ of degree $>2p_a(\Gamma)$, $(E_H)_{|\Gamma}$ is stable (see \cite{EL}).

Since $E_H$ is $\mu-$stable with respect to $2H$, it is also $\mu-$stable with respect to $H$ and this yields the $\mu-$stability of $E_L$ with respect to $L$, because $\pi$ is a double covering (see \cite{huy}, Lemma 3.2.2).
\qed\\

\textbf{Remark.} Throughout the proof the ampleness of $L$ is used only when $C$ is a smooth plane quintic of genus $g=6$ to show that we cannot have equality between slopes. Indeed, if we only assume that $L$ is generated by its global sections and $ L^2\geq 2 $ then $E_L$ is still $\mu-$semistable with respect to $L$ and also $\mu-$stable unless $C$ is a smooth plane quintic of genus $g=6$. 

\section{About abelian surfaces}\label{s6}

In this section we study the same problem when $X$ is an abelian surface over $\mathbb{C}$ and we give the proof of Theorem \ref{abel}.

\begin{pro}\label{specserab}
Let $X$ be an abelian surface over $\mathbb{C}$; then there is no irreducible hyperelliptic curve of genus $g\geq 6$ and no irreducible trigonal curve of genus $g\geq 8$ on $X$.
\end{pro}
\textbf{Proof.} Take $d=2$ or $3$ and suppose that there is a $d-$gonal irreducible curve $C$ of genus $g\geq 2d+2$ on $X$.
Then there is an exact sequence of sheaves on $X$
$$
\xymatrix{
0\ar@{->}[r]&F^*\ar@{->}[r]&H^0(g^1 _d)\otimes \of _X\ar@{->}[r]&\of _C(g^1 _d)\ar@{->}[r]&0}
$$
where $F$ is a vector bundle of rank 2 such that $c_1(F)=C$ and $c_2(F)=d $. Dualising the above exact sequence we get
$$
\xymatrix{
0\ar@{->}[r]&\of _X ^2\ar@{->}[r]&F \ar@{->}[r]&\of _C(K_C-g^1 _d)\ar@{->}[r]&0}
$$
It follows from the assumption on the genus that $c_1(F)^2-4c_2(F)=2g-2-4d>0$, so $F$ is Bogomolov unstable (see \cite{Ray}). Therefore, there exists a line bundle $\of _X(A)$ on $X$ such that $\mu(\of _X(A))>\mu(F)$, i.e. $2A.C>C^2$, and we have an exact sequence
$$
\xymatrix{
0\ar@{->}[r]&\of _X(A)\ar@{->}[r]&F \ar@{->}[r]&\mathcal{I}_Z\otimes\of _X(B)\ar@{->}[r]&0}
$$

with $A+B=C$, $A.B+\deg \mathcal{I}_Z=d$ and $(A-B)^2>0$ (see \cite{Ray}). Hence we can construct the following diagram
$$
\xymatrix{
&&0\ar@{->}[d]&&\\
&&\of _X ^2\ar@{->}[d]^i\ar@{->}[dr]&&\\
0\ar@{->}[r]&\of _X(A)\ar@{->}[r]&F \ar@{->}[r]\ar@{->}[d]&\mathcal{I}_Z\otimes\of _X(B)\ar@{->}[r]&0\\
&&\of _C(K_C-g^1 _d)\ar@{->}[d]&&\\
&&0&&\\}
$$
Since $i$ is an isomorphism outside $C$, $h^0(\mathcal{I}_Z\otimes\of _X(B))>0$ and $B$ is effective. By the Hodge index theorem $A^2B^2\leq (A.B)^2\leq d^2$. Since $K_X=0$, $A^2 $ and $B^2$ are even numbers and $A^2>B^2$ because $2A.C>C^2$, hence we must have $B^2\leq 2$.

If $B^2=2$, then $ d=3 $ and $A^2=4$ and we would have $6-2A.B>0$, so $A.B\leq 2$ in contradiction with  $A^2B^2=8$. Therefore $B^2=0$, which means that $B=kE$ where $E$ is an elliptic curve and $k\geq 1$; on the other hand we know that $0\leq A.B\leq d$. In fact $A.B>0$, otherwise by the Hodge index theorem it would follow $B=0$ against the fact that $h^0(\mathcal{I}_Z\otimes\of _X(B))>0$; hence $1\leq kA.E\leq d$. Since $A.E=1$ would imply that $A$ itself is elliptic, the only possibility is $k=1$ and $A.B>1$. In this case we have $h^0(B)=1$, hence by the snake lemma we have the following diagram
$$
\xymatrix{
&0\ar@{->}[d]&0\ar@{->}[d]&0\ar@{->}[d]&\\
0\ar@{->}[r]&\of _X \ar@{->}[d]^s\ar@{->}[r]&\of _X ^2\ar@{->}[d]\ar@{->}[r]&\of _X \ar@{->}[d]^{\sigma}\ar@{->}[r]&0\\
0\ar@{->}[r]&\of _X(A)\ar@{->}[r]\ar@{->}[d]&F \ar@{->}[r]\ar@{->}[d]&\mathcal{I}_Z\otimes\of _X(B)\ar@{->}[r]\ar@{->}[d]&0\\
0\ar@{->}[r]&\tau \ar@{->}[r]\ar@{->}[d]&\of _C(K_C-g^1 _d)\ar@{->}[d]\ar@{->}[r]&\tau '\ar@{->}[r]\ar@{->}[d]&0\\
&0&0&0&\\}
$$
where $\tau$ and $\tau'$ are two torsion sheaves with support respectively on the zero-locus of $s$ and $\sigma$.
Hence the exactness of the third line implies that $C$ is reducible, against our assumptions. \qed\\

\textbf{Proof of Theorem \ref{abel}.} 
Since $L$ is generated by its global sections such that $L^2\geq 14$, the general member of $|L|$ is a smooth irreducible curve of genus $g\geq 8$. Hence, given a non-zero $\alpha\in \mathrm{Pic}^0(X)$, we can find $C\in |L\otimes\alpha ^{-1}|$ smooth irreducible of genus $g\geq 8$. The $\mu-$stability of $E_L$ with respect to $L$ is equivalent to the $\mu-$stability of $E_L$ with respect to $C$. Since we have $H^0(\alpha)=H^1(\alpha)=0$, it follows from Proposition \ref{pdiag} that $(E_L)_{|C}\cong E_{(L_{|C})}$. Moreover, $L_{|C}\cong K_C\otimes \alpha _{|C}$, so  by Theorem \ref{p:pa} $ E_L $ is $ \mu- $stable with respect to $ C$ if $ c(C)\geq 2 $.
By the hypothesis on the genus of $C$ and by Proposition \ref{specserab} the cases $ c(C)=0,1 $ cannot occur, so there is nothing more to prove.\qed\\

\textbf{Remark.} In the case $ g(C)\leq 7 $ the same proof shows the $\mu-$stability of $E_L$ if $c(C)\geq 2$. Moreover, it is possible to show that $ E_L $ is $\mu-$stable with respect to $L$ also if either $C$ is a smooth plane quintic of genus $g=6$ or if $C$ is a trigonal curve of genus $g=4$.
\section*{Acknowledgements}
I am very grateful to my advisor Prof. Arnaud Beauville for the help he gave me throughout this year and for patiently reading all the drafts of this paper.

\end{document}